\documentclass[12pt]{article}
\usepackage{amsmath}
\usepackage{amssymb}

\newtheorem{thm}{Theorem}[section]
\newtheorem{lem}{Lemma}[section]

\numberwithin{equation}{section}

\def\bbR{{\mathbb R}}
\def\bbZ{{\mathbb Z}}
\def\bbN{{\mathbb N}}
\def\bbQ{{\mathbb Q}}
\def\bbE{{\bf E}}
\def\bbP{{\bf P}}

\def\eps{\varepsilon}
\let\phi=\varphi
\def\qed{\hfill $\square$}

\begin{document}

\title{The shape theorem for the frog
model\thanks{The authors are thankful to CAPES/PICD,
CNPq (300226/97--7) and FAPESP (97/12826--6) for financial
support.}}

\author{O.S.M.~Alves \and F.P.~Machado \and S.Yu.~Popov}

\maketitle {\footnotesize \noindent Department of Statistics,
Institute of Mathematics and Statistics, University of S\~ao
Paulo, Rua do Mat\~ao 1010, CEP 05508--900, S\~ao Paulo SP,
Brazil

\noindent e-mails: \texttt{oswaldo@ime.usp.br, fmachado@ime.usp.br,
popov@ime.usp.br} }

\begin{abstract}
In this work we prove a shape theorem for a growing set of
Simple Random Walks (SRWs) on $ \bbZ^d$, known as frog model.
The dynamics of this process is described as follows: There are
active particles, which perform independent SRWs, and sleeping
particles, which do not move. When a sleeping particle is hit by an
active particle, it becomes active too.
At time~$0$ all particles are sleeping, except
for that placed at the origin. We prove that the set of the
original positions of all the active particles, rescaled by the
elapsed time, converges to some compact convex set. In some specific cases
we are able to identify (at least partially) this set.
\\[.3cm] {\bf Keywords:} frog model, shape theorem,
subadditive ergodic theorem
\end{abstract}

\section{Introduction and results}
\label{intro}
We study a discrete time particle system on $ \bbZ^d $ named frog
model. In this model particles, thought of as frogs, move as
independent simple random walks (SRWs) on $ \bbZ^d $.
At time zero there is one particle at each site of the
lattice and all the particles are sleeping except for the one at
the origin. The only awakened particle starts to perform a SRW.
  From then on when an awakened particle jumps on the top of a
sleeping particle, the latter wakes up and starts jumping
independently, also performing a SRW. The number of awakened
particles grows to infinity as active particles jump on sites that
have not been visited before, awakening the particles that are
sitting there. Let us underline that the active particles do not
interact with each other and there is no
``one-particle-per-site'' rule.

This model is a modification of a model for information
spreading that the authors have learned from K.~Ravishankar.
The idea is that every moving
particle has some information and it shares that information with
a sleeping particle at the time the former jumps on the back of
the latter. Particles that have the information move freely
helping in the process of spreading information.
The model that we deal with in this paper
is a discrete-time version of that proposed
by R.~Durrett, who also suggested the term ``frog model''.

In~\cite{P} this model was studied from the point of view of
transience and recurrence in the case when the initial
configuration of sleeping particles is random and has decaying
density.

We now define the process in a more formal way. Let $ \{(S^x_n)_{n
\in \bbN}, {x \in \bbZ^d}\} $ be independent SRWs such that $
S^x_0 = x $ for all $ x \in \bbZ^d $. For the sake of cleanness let $
S_n := S^0_n $. These sequences of random variables give the
trajectory of the particle seated originally at site $ x $, starting
to move at the time it wakes up. Let
\begin{equation}
\label{tht}
t(x,z) = \min \{n: S^x_n =z \},
\end{equation}
remembering that if $d > 2$ then $\bbP[t(x,z) = \infty] > 0 $.

We define now
\begin{equation}
\label{Tht}
T(x,z) = \inf \Big\{\sum_{i=1}^k t(x_{i-1}, x_i): x_0 = x,
\ldots ,x_k=z \mbox{ for some } k \Big\},
\end{equation}
 the passage time from $x$ to $z$ for the frog model.
By now it should
be clear that if the process starts from just one active particle
sitting at site~$x$, in the sense that that particle is the only
active one, $ T(x,z) $ is the time it takes to have the
particle sitting at site~$z$ to be awakened. Note that the particle
which awakens~$z$ need not be that from~$x$.

Now, let $Z_y^x (n)$ be the location (at time $n$)
 of the particle that started from site~$y$
 in the process in which the only active particle
at time zero was at site~$x$. Formally, we have
$ Z_x^x (n) = S_n^x $, and
\[
 Z_y^x (n) = \left\{
  \begin{array}{ll}
  y, & \mbox{if } T(x,y) \ge n,\\
  S_{n-T(x,y)}^y, & \mbox{if }  T(x,y) < n.
  \end{array}
  \right.
\]
Since every random variable of the form $Z_y^x (n)$ is constructed using
the random variables $\{(S^x_n)_{n\in \bbN}, {x \in \bbZ^d}\}$
and the related random variables $\{ T(x,y): x,y \in \bbZ^d\}$, this defines a
coupling of processes $\{Z^x, x\in\bbZ^d\}$.
We point out that a particle starting from a
given site, as soon as it wakes up (becomes active), executes the same
random walk in all the processes. The only difference for that particle is
in the time when it wakes up.

With the help of these variables we define the sites whose
originally sleeping particles have been awakened by time~$n$,
by some originally awakened particle from the set $ \mathsf{B} $.
Namely
\[
 \xi^{\mathsf{B}}_n = \{ y \in \bbZ^d : T(x,y) \le n \mbox{ for some }
x \in \mathsf{B} \}.
\]
We are mostly concerned with $ \xi_n := \xi^0_n $
(note that the process~$\xi_n$ itself is {\it not\/} Markovian).
Moreover, we define
\[
 {\bar \xi}_n = \{ x + ({-1/2}, {1/2}]^d: x \in \xi_n\}\subset\bbR^d.
\]

It is a basic fact that the displacement of a single SRW at time~$n$
is roughly $\sqrt n$. Here we prove that in the frog model the
particles may ``help'' each other in order to make the
boundary of ${\bar \xi}_n$ grow linearly with time. Namely,
the main result of this paper is the following

\begin{thm}
\label{maintheo} There is a non-empty convex set $ {\mathsf A} \subset
\bbR^d $ such that, for any $ 0 < \eps < 1 $
\[
(1-\eps){\mathsf A} \subset
\frac{{\bar \xi}_n}{n} \subset (1+\eps){\mathsf A}
 \]
 for all $ n $ large enough almost surely.
\end{thm}

Note that, although Theorem~\ref{maintheo} establishes the existence of the
asymptotic shape~${\sf A}$, it is difficult to say something
definite about this shape, except, of course, that~${\sf A}$ is
symmetric and that ${\sf A}\subset {\cal D}$, where
\[
{\cal D} = \{x=(x^{(1)},\ldots,x^{(d)})\in \bbR^d :
     |x^{(1)}|+\cdots +|x^{(d)}| \leq 1\}.
\]
Also, note that if the initial configuration is augmented (i.e.\
some new particles are added), then the asymptotic shape (when it
exists) augments as well. We are going to show that if the
initial configuration is rich enough, then the limiting
shape~${\sf A}$ may contain some pieces of the boundary of~${\cal
D}$ (a ``flat edge'' result) or even coincide with~${\cal D}$ (a
``full diamond'' result).

To formulate these results, we need some additional notation. For $1\leq
i<j\leq d$ let
\[
 \Lambda_{ij} = \{x=(x^{(1)},\ldots,x^{(d)})\in \bbR^d :
      x^{(k)}=0 \mbox{ for } k\neq i,j\},
\]
and for $0<\beta<1/2$ let
\[
\Theta^\beta_{ij} = \{x=(x^{(1)},\ldots,x^{(d)}) \in \Lambda_{ij} :
 |x^{(i)}|+|x^{(j)}|=1, \min\{|x^{(i)}|,|x^{(j)}|\}\geq \beta\}.
\]
Define $\Theta^\beta$ to be the convex hull of
$(\Theta^\beta_{ij})_{1\leq i<j\leq d}$. Denote by~${\sf A}_m$ the
asymptotic shape in the frog model when the initial configuration
is such that any site~$x\in \bbZ^d$ contains exactly~$m$
particles. The existence of~${\sf A}_m$ for arbitrary~$m$ can be
derived in just the same way as in the case $m=1$
(Theorem~\ref{maintheo}). Then, for a positive integer-valued
random variable~$\eta$ denote by~${\sf A}_\eta$ the asymptotic
shape (if exists) for the frog model with the initial
configuration constructed in the following way: Into every site we
put a random number of particles independently of other sites, and
the distribution of this random number is that of~$\eta$.

\begin{thm}
\label{flat_edge} If $ m $ is large enough, then there exists
$0<\beta<1/2$ such that $\Theta^\beta \subset {\sf A}_m$.
\end{thm}

\begin{thm}
\label{full_diamond} Suppose that for some positive $\delta<d$ and
for all~$n$ large enough we have
\begin{equation}
\label{fd}
 \bbP[\eta\geq (2d)^n] \geq n^{-\delta},
\end{equation}
and $\eta\geq 1$ a.s. Then ${\sf A}_\eta = {\cal D}$.
\end{thm}

The paper is organized in the following way. Section~\ref{basics}
contains some well known results about large deviations and
SRW on $ \bbZ^d $. We need
these results to verify the hypotheses of Liggett's subadditive
ergodic theorem. These hypotheses are verified in Section~\ref{subadd} and
the proof of the shape theorem is given in Section~\ref{AS}. Besides
that, in Section~\ref{AS} we prove the ``flat edge'' and
the ``full diamond'' results.

\section{Basic facts}
\label{basics}
Along this section we state basic facts about large deviations and
random walks which we need to prove our results. A couple of them
are followed by their proofs just because we have not been able to
find them in the bibliography. As usual, $ C, C_1, C_2, \ldots $
stand for positive finite constants. For what follows we use these
constants freely. Also, $\lfloor x \rfloor$ stands for the largest
integer which is less than or equal to~$x$, while
$\lceil x \rceil$ is the smallest
integer which is greater than or equal to~$x$.

\begin{lem}
\label{xex}
If $ X $ is a random variable assuming positive integer values
such that $ X \le a $ almost surely and $ \bbE X \ge b > 0 $ then
\begin{equation}
\label{basicprob}
\bbP\Big[X\ge \frac{b}{2}\Big] \ge \frac{b}{2 a}.
\end{equation}
\end{lem}

\noindent {\it Proof.} We have
\begin{eqnarray*}
 b \le \bbE X &=& \sum_{i<b/2} i \bbP[X=i]
+ \sum_{i\geq b/2} i \bbP[X=i] \\
 &\le & \frac{b}{2} +  a \bbP\Big[X \ge
\frac{b}{2}\Big],
\end{eqnarray*}
and~(\ref{basicprob}) follows.
\qed

\subsection{On large deviations}
In this subsection we state two large deviations results. The
first one (Lemma \ref{Shiryaev}) is a simple
application of the exponential Chebyshev inequality
to sums of independent
Bernoulli random variables. The second one
(Lemma~\ref{Nagaev}) is useful when one has positive integer
random variables but cannot guarantee the existence of their
moment generating functions. That condition is weakened and
substituted by a sub exponential estimate for the tails of their
probability distributions.

\begin{lem} [\cite{Shiryaevp},  p.\ 68.]
\label{Shiryaev} Let $ \{X_i, i \ge 1\} $ be i.i.d.\ random variables with
$\bbP[X_i=1] = p$ and $\bbP[X_i=0] = 1-p$. Then for any
$0<p<a<1$ and for any $N \geq 1$ we have
\begin{equation}
\label{shi} \bbP\Big[\frac{1}{N}\sum_{i=1}^N X_i \geq a\Big]
\leq \exp\{-N H(a,p)\},
\end{equation}
where
\begin{equation}
\label{H} H(a,p) = a \log\frac{a}{p} + (1-a)\log\frac{1-a}{1-p} >
0.
\end{equation}
\end{lem}
Next large deviation result is an immediate consequence
of Theorem~1.1, p.~748 of~\cite{Nagaev}.

\begin{lem}
\label{Nagaev} Let $ \{X_i, i \ge 1\} $ i.i.d.\ positive integer-valued
random variables such that there are $ C_1 > 0 $ and $ 0 < \alpha < 1 $
such that
\begin{equation}
\label{subexponential}
\bbP[X_i \geq n] \le C_1 \exp\{-n^{\alpha}\}.
\end{equation}
Then there exist $ a > 0, \ 0 < \beta < 1 $ and $ C_2 > 0 $ such
that for all~$n$
\[
\bbP\Big[ \sum_{i=1}^n X_i  \geq a n \Big] \le C_2 \exp \{-n^{\beta}\} .
\]
\end{lem}

\subsection{On simple random walks}

The following results for $d-$dimensional SRW are found
in~\cite{Hughes} and~\cite{Lawler}. Let $ {\mathsf p}_n(x) =
\bbP[S_n = x] $ and $ \|x\| $ be the Euclidean norm.
   From Section~\ref{subadd} onwards we
also work with the norm $ \|x\|_1 = \|(x^{(1)}, \ldots, x^{(d)})
\|_1 = \sum_{i=1}^d |x^{(i)}|$.

\begin{thm} [\cite{Lawler},  p.\ 14, 30.]
\label{pnx}
\begin{equation}
\label{pnx1} {\mathsf p}_n(x) = \frac{2}{n^{d/2}} \Big( \frac{d}{2
\pi} \Big)^{d/2} \exp\Big\{\frac{-d\|x\|^2}{2n}\Big\}+e_n(x),
\end{equation}
where $ |e_n(x)| \le Cn^{-(d+2)/2} $ for some $ C $
and for all $ a < d $
\[
 \lim_{\|x\| \to \infty} \|x\|^{a} \sum_{n=1}^{\infty} |e_n(x)|
= 0 .
\]
\end{thm}

\begin{thm} [\cite{Lawler},  p.\ 29.]
\label{psupsrw}
There is $ C $ such that for all $ n, t $
\begin{equation}
\bbP [ \sup_{0 \le i \le n}\|S_i\| \ge t n^{1/2}] \le Ce^{-t}.
\end{equation}
\end{thm}

Let
\[
 {\mathsf R}^{\mathsf B}_n = \{ S^{\mathsf B}_i : 0 \le i \le n \} =
\{ y \in \bbZ^d : t(x,y) \le n
\mbox{ for some } x \in {\mathsf B} \}
\]
be the set of distinct sites visited by
 the family of SRWs, starting from the
set of sites $ \mathsf B $, up to time~$ n $. Some authors refer to
${\mathsf R}^0_n $ as the range of SRW. As usual,
$ |{\mathsf R}^{\mathsf B}_n| $ stands for the
cardinality of $ {\mathsf R}^{\mathsf B}_n $. A useful basic fact is
that $ |{\mathsf R}^{\mathsf B}_n| \le (n+1) {|{\mathsf B}|} $.

\begin{thm} [\cite{Hughes},  p.\ 333, 338.]
\label{esprange}
 \begin{itemize}
\item If $ d = 2 $ then there is $ a_2 >0 $ such that
\begin{equation}
 \lim_{n \to \infty} \frac{\bbE |{\mathsf R}^0_n| }{{n / \log n}} = a_2.
 \end{equation}
\item If $ d \ge 3 $ then there is $ a_3 := a_3(d) > 0 $ such that
\begin{equation}
\lim_{n \to \infty} \frac{\bbE |{\mathsf R}^0_n|}{n} = a_3.
\end{equation}
\end{itemize}
\end{thm}

Let $ G_n(x) := \sum_{j=0}^n {\mathsf p}_j(x) $ be the mean
number of visits to site $ x $ up to time~$ n $ and $ G(x) =
G_{\infty}(x) $. These are the well known {\it Green's functions}.
Let $ {\mathsf q}(n,x) = \bbP[t(0,x) \le n] $. From
Theorem~\ref{pnx} we get the following result:

\begin{thm}
\label{pofvisit}
\begin{itemize}
\item If $ d = 2, x \not= 0 $ and $  n \geq {\|x\|}^2 $,
then there exists $ C_2 > 0 $ such that
\begin{equation}
{\mathsf q}( n, x ) \ge \frac{C_2}{\log \|x\|} .
\end{equation}
\item Suppose that $ d \ge 3,  x \not= 0 $ and $ n \ge \|x\|^2 $.
Then there exists $ C_3 = C_3(d) > 0 $ such that
\begin{equation}
{\mathsf q}( n, x) \ge \frac{C_3}{\| x \|^{d-2}} .
\end{equation}
\end{itemize}
\end{thm}

\noindent {\it Proof.} Suppose without loss
of generality that $\|x\|^2 \leq n \leq \|x\|^2+1$. Observe that
\begin{eqnarray*}
G_n(x) &=& \sum_{j=0}^n {\mathsf p}_j(x) = \sum_{j=0}^n
\sum_{k=0}^j {\mathsf p}_k(0) \bbP[t(0,x)=j-k] \\
&=& \sum_{k=0}^n {\mathsf p}_k(0) {\mathsf q}(n-k,x)
\le {\mathsf q}(n,x)G_n(0).
\end{eqnarray*}
So
\[
 {\mathsf q}(n,x) \ge \frac{G_n(x)}{G_n(0)} \ge
\left\{
\begin{array}{ll}
\displaystyle
\frac{\sum_{j={\lfloor n/2 \rfloor}}^n
\mathsf{p}_j(x)}{\sum_{j=0}^n{\mathsf p}_j(0)}, & d=2, \\
\vphantom{\sum^{\int^N}}
(G(0))^{-1} \sum_{j={\lfloor n/2 \rfloor}}^n
\mathsf{p}_j(x), & d \geq 3.
\end{array}
\right.
\]
 Using~(\ref{pnx1}), after some elementary computations
we finish the proof. \qed

\section{Subadditive ergodic theorem}
\label{subadd}
The basic tools for proving shape theorems are the subadditive
ergodic theorems. Next we state Liggett's subadditive ergodic
theorem (cf., for example, \cite{D}).
 In the sequence we apply it to the random variables $T(\cdot,\cdot)$.

\begin{thm}
\label{subaet} Suppose that $ \{ Y(m,n) \} $ is a collection of
positive random variables indexed by integers satisfying $ 0 \le m
< n $ such that

\begin{itemize}
\item $ Y(0,n) \le Y(0,m) + Y(m,n)$
   for all $ 0 \le m < n $ (subadditivity);

\item The joint distribution of $ \{Y(m+1, m+k+1), k \ge 1\} $
is the same as that of $ \{ Y(m, m+k), k \ge 1\} $ for each $ m
\ge 0 $;

\item For each $ k \ge 1 $ the sequence of random variables
$ \{Y(nk, (n+1)k), n \ge 1\} $ is a stationary
ergodic process;

\item $ \bbE Y(0,1) < \infty $.
\end{itemize}
Then
\[
 \lim_{n \to \infty} \frac{Y(0,n)}{n} \to \gamma \qquad \mbox{a.s.},
\]
where
\[
 \gamma = \inf_{n \ge 0} \frac{\bbE
Y(0,n)}{n}.
\]
\end{thm}

In order to verify the hypotheses of Liggett's subadditive
ergodic theorem, for each fixed $ x \in \bbZ^d $ let us
consider $ Y(m,n) = T(mx,nx) $.

First of all observe that the set of variables $ \{T(x,y):x,y \in \bbZ^d \} $
defined in Section~\ref{intro} is {\it subadditive\/} in the sense that
\begin{equation}
\label{subadit}
T(x,z) \le T(x,y) + T(y,z) \mbox{ for all } x, y, z
\in \bbZ^d.
\end{equation}
Instead of proving this fact in a formal way, we prefer to
give a verbal explanation. If site~$z$
is reached before site~$y$, there is nothing to prove. If that
does not happen, observe that the process
departuring from only site~$y$ awakened (the one which gives the
passage time $T(y,z)$) is coupled with the original
process, for which one might have
other particles awakened at time $T(x,y)$ besides that from~$y$.
Thus the remaining time to reach site~$z$ for the original
process can be at most $T(y,z)$. This takes care of the first hypothesis.

For a fixed $ x \in \bbZ^d $ and $ k \in \bbN $, the sequence
$ \{T((n-1)kx, nkx) , n \ge 1 \} $ is {\it stationary} by definition.
{\it Ergodicity\/} holds because the sequence of random variables
$ \{T((n-1)kx, nkx) ,  n \ge 1 \} $ is strongly mixing. That can be checked
easily because the
events $ \{ T(n_1 kx, (n_1 +1)kx) = a \}  $ and
$ \{ T(n_2 kx, (n_2 +1)kx) = b \} $ are independent
provided that $ a+b < \|(n_1-n_2)kx\|_1 $.

It is simple to see that the fourth hypothesis holds when $ d =
1$. To see that remember that for $ \tau = $ {\it the first return
to the origin of a SRW}, we can assure that
$ \bbP [ \tau > t ] \le C t^{-1/2} $.
Besides that, in a random time with exponential tail we will
have at least three awakened particles jumping independently in the
frog model. Combining these two facts we have that
$ \bbE T(0,1) < \infty$. So, for $d=1$ one gets $T(0,n)/n \to \gamma$
a.s., and consequently we have the proof of the shape
theorem with $\mathsf{A} = [-\gamma^{-1},\gamma^{-1}]$ in dimension~$1$.

For higher dimensions we need a more
powerful machinery in order to check the fourth hypothesis.
This comes in the next result.

\begin{thm}
\label{expdecay}
For all $ d \ge 2 $ and $ x_0 \in \bbZ^d $ there exist positive
finite constants $ C=C(x_0) $ and $ \gamma $ such that
\[
 \bbP [ T(0,x_0) \ge n ] \le  C \exp\{- n^{\gamma}\}.
\]
\end{thm}

\noindent {\it Proof.} For technical reasons we treat the case $ d
= 2 $ separately. Let $d=2$. First pick $ n \ge \|x_0\|^2 $.
Remember that
\[
{\mathsf R}^0_n = \{ x \in \bbZ^d : t(0,x) \le n \}
\]
 is the range of the SRW of the first awakened particle,
up to time $ n $. By Theorem~\ref{esprange} we have that for all~$k$
large enough it is true that $\bbE|{\mathsf R}^0_k| \ge {2 C_1 k / \log k}$.
 Since  $ |{\mathsf R}^0_k| \le k + 1 $, by
using Lemma~\ref{xex} we obtain
\begin{equation}
\label{p2range}
\bbP\Big[|{\mathsf R}^0_k| \ge \frac{C_1 k}{\log k}\Big]
  \ge \frac{C_2}{\log k}.
\end{equation}

Divide the time interval $ [0,n] $ into (roughly speaking) $
n^{1/2} $ disjoint intervals of size $ n^{1/2} $.
Next we keep track of the displacement of the
original particle over that time interval. For each subinterval of
size $ n^{1/2} $, the cardinality of the corresponding
 subrange does not depend
on the cardinalities of other subranges. Therefore, for the event
\[
 A := \Big\{ | \xi_n | \ge \frac{2 C_1 n^{1/2}}{\log n} \Big\} \supset
\Big\{ | {\mathsf R}^0_n | \ge \frac{2 C_1 n^{1/2}}{\log n} \Big\}
\]
by~(\ref{p2range}) it holds that
\begin{equation}
\label{pprobA}
 \bbP[A] \ge 1-\Big(1- \frac{2 C_2}{\log n}\Big)^{n^{1/2}}.
\end{equation}
Let us now consider the event
\[
 B = B(n, \eps) = \{\sup_{0 \le i \le n} \|S^0_i \| <
n^{{1/2}+\eps} \}
\]
where $ 0 < \eps < 1$.
Observe that by Theorem~\ref{psupsrw} there is  $ C_3 $ such that
\begin{equation}
\label{pprobBc}
 \bbP[B^c] \le C_3 \exp\{-n^{\eps}\}.
\end{equation}
Note that
\begin{equation}
\label{ur*}
\bbP[T(0,x_0) > n + 4 n^{1 + 2 \eps} ] \le
\bbP[T(0,x_0) > n + 4 n^{1 + 2 \eps} | A \cap B ]
+ \bbP[A^c] + \bbP[B^c].
\end{equation}
Considering now all the particles awakened by the
first particle, with the help of Theorem~\ref{pofvisit} one gets
\begin{equation}
\label{ur**}
 \bbP[T(0,x_0) > n + 4 n^{1 + 2 \eps} | A \cap B ] \le
\Big(1-\frac{C_3}{\log n}\Big)^{2 C_1 n^{1/2}(\log n)^{-1}}
\end{equation}
 for any fixed $ 0 < \eps < 1$, so the result follows for $d= 2$
  from (\ref{pprobA})--(\ref{ur**}).


Now we treat higher dimension cases. For each dimension $ d \ge 3
$ fixed, the proof needs $ \lfloor d/2 \rfloor + 1 $ steps. To do
that in a general fashion, we separate the proof in four parts
named {\it first step, second step, general reasoning\/} and {\it
denouement}. As the case $ d = 3 $ is simpler, we are able to
finish its proof right after the {\it first step}. For $ d = 4, 5
$ we skip the part {\it general reasoning}, going directly to the
part {\it denouement}. For higher dimensions all parts are needed.

\medskip \noindent {\it First step:}

\noindent
Pick $ n \ge \|x_0\|^2 $. By Theorem~\ref{esprange} and Lemma~\ref{xex}
it follows that for some $r_1>0$
\begin{equation}
\label{pranfg}
\bbP[|{\mathsf R}^0_k| \ge r_1 k] \ge C_1.
\end{equation}
Divide the time interval $ [0,n] $ into (roughly
speaking) $ n^{\eps} $ disjoint intervals of size $
n^{1-\eps} $ for some $ 0 < \eps < 1 $ to be defined later.
Next we keep track of the range of the
original particle over that time interval. For each subinterval of
size $ n^{1-\eps}$, the cardinality of the
corresponding subrange does not depend
on the cardinalities of other subranges. Let
\[
A_1 := A_1(n, \eps) := \{ | \xi_n | \ge r_1 n^{1- \eps}\}
\supset \{ | {\mathsf R}^0_n | \ge r_1 n^{1- \eps}\} ;
\]
 by~(\ref{pranfg}) it holds that
\begin{equation}
\label{paum}
\bbP[A_1] \ge 1-(1-C_1)^{n^\eps}.
\end{equation}
Again consider the event
\[
 B = B(n, \eps)= \{\sup_{0 \le i \le n} \|S^0_i \| <
n^{{1/2}+\eps} \}
\]
and observe that by Theorem~\ref{psupsrw} there is  $ C_2 $ such that
\begin{equation}
\label{probBc}
 \bbP[B^c] \le C_2 \exp\{-n^{\eps}\}.
\end{equation}
 For $ d = 3 $ we are done, since by
Theorem~\ref{pofvisit}, analogously to~(\ref{ur*})--(\ref{ur**}) we have
\[
 \bbP[T(0,x_0) > n + 4 n^{1 + 2 \eps} \mid A_1 \cap B ] \le
\Big(1-\frac{C_3}{n^{{1/2} + \eps}}\Big)^{r_1 n^{1-\eps}}.
\]
By choosing $ \eps < {1}/{4}$ and using (\ref{paum})--(\ref{probBc}),
 the result follows for $d=3$.

\medskip
\noindent
{\it Second step:}

\noindent
Denote
\[
 \mathsf{D}_{i, \eps} := \{x\in\bbZ^d : \|x\| \le i n^{{1/2}+ \eps}\}.
\]
Suppose that events $ A_1 $ and $ B $ happen.
Suppose also that all the particles awakened by the original
particle until time~$n$ are allowed to start moving {\it exactly\/}
at the moment~$n$; clearly, such a procedure can only increase
the hitting time of~$x_0$ for the frog model. In this case,
there are {\it at least} $r_1 n^{1-\eps}$ active particles at
time $ n $ and all of them are inside the ball
$ \mathsf{D}_{1, \eps} $. Choose $ r_1 n^{1-\eps} $ of these
particles and call this set $ G_1 $. Split $ G_1 $ into $ r_1
n^\eps $ groups of size $ n^{1-2\eps} $. Call these groups
$ G_1^1, G_1^2, \ldots , G_1^{r_1 n^\eps}$. For each~$y$ in
the ring $ \mathsf{D}_{2, \eps} \setminus \mathsf{D}_{1, \eps} $
let $ \zeta^{(2)}_i(y) $ be the indicator function of the following event
\[
 \{ \mbox{there exists }
x \in G^1_1 \mbox{ such that } t(x,y) \le 9 n ^{1 + 2 \eps} \}.
\]
By Theorem~\ref{pofvisit}, and using the fact that for
$ x \in {\mathsf D}_{1, \eps} $,
$ y \in {\mathsf D}_{2, \eps} \setminus {\mathsf D}_{1, \eps} $
the distance $ \|x-y\| $ is less than or equal to
$ 3n^{1/2 + \eps}$,
we have
\begin{eqnarray}
\label{expislarge}
\bbE(\zeta^{(2)}_i(y) | A_1 \cap B )
&=& \bbP[\zeta^{(2)}_i(y) = 1 | A_1 \cap B]
\nonumber \\
&\ge& 1 - \prod_{x \in G^i_1} (1-{\mathsf q}(9n^{1+2 \eps},y-x))
\nonumber \\
&\ge& 1-\Big(1-
\frac{C_3}{3^{d-2}n^{({1/2}+\eps)(d-2)}}\Big)^{n^{1-2\eps}}
\nonumber \\
&\ge& \frac{C_4}{n^{{d/2}+d \eps - 2}}
\end{eqnarray}
(here we used the fact that $d/2+d\eps > 2$ for $d\geq 4$).
 Let
\[
 \zeta^{(2)}_i = \sum_{y \in \mathsf{D}_{2, \eps}
\setminus \mathsf{D}_{1, \eps}} \zeta^{(2)}_i(y).
\]
Since $|\{ x \in \bbZ^d : x \in \mathsf{D}_{2, \eps}
\setminus \mathsf{D}_{1, \eps}
\}| \ge C_5 n^{{d/2}+d\eps} $, we have that there exists~$r_2>0$
such that
\begin{equation}
\label{ur'*}
 \bbE(\zeta^{(2)}_i | A_1 \cap B) \ge 2 r_2 n^2.
\end{equation}
Clearly, it is true that
\begin{equation}
\label{ur'**}
 \zeta^{(2)}_i \le n^{1-2\eps}\times 4n^{1+2\eps}= 4 n^2 .
\end{equation}
Therefore, using Lemma~\ref{xex},
we obtain from (\ref{ur'*})--(\ref{ur'**}) that
\begin{equation}
\label{jkl}
 \bbP[\zeta^{(2)}_i \ge r_2 n^2|A_1 \cap B ] \ge C_6 > 0.
\end{equation}
For $ n_2 := n + 9n^{1+2 \eps} $ let
\[
A_2 = A_2(n, \eps) = \{ | \xi_{n_2} \cap
(\mathsf{D}_{2, \eps} \setminus \mathsf{D}_{1, \eps}) |
\ge r_2 n^2 \} .
\]
 Repeating the above argument for the groups $G_1^2,\ldots,
G_1^{r_1 n^{\eps}}$, we
obtain from (\ref{jkl}) that
\begin{equation}
\label{d45}
 \bbP[A^c_2 | A_1 \cap B] \le (1-C_6)^{r_1 n^{\eps}}.
\end{equation}
 For $ d=4 $ or $5$, one should go directly to {\it
denouement}.

\medskip
\noindent
{\it General reasoning:}

\noindent
Consider the sequence of times
\[
 n_1 := n , n_2 := n + 9n^{1 + 2 \eps}, \ldots,
n_k := n + n^{1+2 \eps} \sum_{j=2}^k (2j-1)^2,
\]
$k\leq \lfloor d/2 \rfloor$.
With that sequence, we define the events
\[
 A_k = \{ | \xi_{n_k} \cap
( \mathsf{D}_{k, \eps} \setminus\mathsf{D}_{k-1, \eps} )
| \ge r_k n^k \},
\]
where the constants~$r_k$, $k>2$, will be defined later,
and the random sets
\[
{\tilde G}_k = \{ x \in \mathsf{D}_{k, \eps}
\setminus\mathsf{D}_{k-1, \eps}: T(0,x) < n_k \}.
\]
for $ k \in \{1, \ldots, \lfloor d/2 \rfloor \}$.
We claim that for $2\leq k\leq \lfloor d/2 \rfloor -1$
\begin{equation}
\label{claim}
\bbP[A_{k+1} | A_k] \ge 1 - \exp\{-Cn^{2\eps}\}.
\end{equation}
 To see why the claim is correct, pick $ r_k n^k $ sites
of $ {\tilde G}_k $ (which are inside of the set
 $ {\mathsf D}_{k, \eps} \setminus
{\mathsf D}_{k-1, \eps} $) at time $ n_k $ and divide them into
$ r_k n^{2 \eps} $ groups of size $ n^{k - 2 \eps} $. Name
these sets $ G_k^1, \ldots, G_k^{r_k n^{2 \eps}} $. Name their
union $ G_k \subset {\tilde G}_k $. As before, we suppose that the
particles that were originally at the sites of set $ G_k $ begin
to move only at time $ n_k $. For each $ y $ in the ring $
{\mathsf D}_{k+1, \eps} \setminus {\mathsf D}_{k, \eps} $
let $ \zeta_i^{(k+1)}(y) $ be the indicator function of the event
\[
\{ \mbox{there exists } x \in G_{k}^1
\mbox{ such that } t(x,y) \le n_{k+1} - n_k \}.
\]
Note that the quantities $n_k$ were defined in such a way
that if $x\in \mathsf{D}_{k, \eps}$, $y\in \mathsf{D}_{k+1, \eps}$,
then $\|x-y\|\leq n_{k+1}-n_k$.
So, by Theorem~\ref{pofvisit}, analogously
to~(\ref{expislarge}), we have
\[
\bbE (\zeta_i^{(k+1)}(y)|A_k) \ge \frac{C_7}{n^{{d/2}+d \eps - (k+1)}}
\]
(note that ${d/2}+d \eps > k+1$ for $k\leq \lfloor d/2 \rfloor -1$).
Let
\[
 \zeta_i^{(k+1)} = \sum_{y \in {\mathsf D}_{k+1}
\setminus {\mathsf D}_{k}} \zeta_i^{k+1}(y).
\]
 Analogously, it follows that there exists $r_{k+1}>0$ such that
\[
 \bbE(\zeta_i^{(k+1)}| A_k ) \ge 2r_{k+1} n^{k+1}
\]
and, clearly,
\[
 \zeta_i^{(k+1)} \le  n^{k - 2 \eps} \times
(2k+1)^2 n^{1+2 \eps}= (2k+1)^2 n^{k+1}.
\]
So, by Lemma~\ref{xex}, there is $ C_8 $ such that
\[
 \bbP[\zeta_i^{(k+1)} \ge r_{k+1} n^{k+1}| A_k ] \ge C_8 > 0 ,
\]
so, considering now all the $r_kn^{2\eps}$ groups, one gets
\[
\bbP [ A^c_{k+1} | A_k ] \le (1-C_8)^{r_k n^{2\eps}},
\]
which by its turn is equivalent to~(\ref{claim}). Now
by~(\ref{paum})--(\ref{probBc}), (\ref{d45})--(\ref{claim})
and using the following inequality
\[
\bbP[A_{\lfloor d/2 \rfloor}] \ge
\bbP[A_{\lfloor d/2 \rfloor} |A_{\lfloor d/2 \rfloor - 1}]
\cdots
\bbP[A_{3} |A_{2}]
\bbP[A_{2} |A_{1}\cap B] \bbP[A_1\cap B]
\]
it follows that
\begin{equation}
\label{condback}
\bbP[A_{\lfloor d/2 \rfloor}] \ge 1- C \exp\{-n^{\gamma_1}\}
\end{equation}
for some $\gamma_1>0$.

\medskip
\noindent
{\it Denouement:}

\noindent
Denote
\[
 I = \sum_{i=2}^{\lfloor d/2 \rfloor} (2i-1)^2.
\]
The idea is to consider the particles
in $ {\tilde G}_{\lfloor d/2 \rfloor} $ at
the moment $ n_{\lfloor d/2 \rfloor} = n + In^{1+2\eps} $
and wait until the moment $ n_{\lfloor d/2 \rfloor} +
d^2n^{1+2 \eps} $ in order to have a overwhelming probability
for them to reach site $ x_0$.

Let
\[
 H := \{ \mbox{no hitting at } x_0 \mbox{ over time interval }
( n_{\lfloor d/2 \rfloor},  n_{\lfloor d/2 \rfloor} + d^2n^{1+2 \eps}]\}.
\]
The number of particles in~$\tilde G_{\lfloor d/2 \rfloor}$
is at least $r_{\lfloor d/2 \rfloor} n^{\lfloor d/2 \rfloor}$ and
they are all at the distance of order $n^{1/2+\eps}$ from~$x_0$, so
 by using Theorem~\ref{pofvisit}, we obtain
\begin{eqnarray*}
 \bbP[ T(0,x_0) > n + (I+d^2)n^{1 + 2 \eps}| A_{\lfloor d/2 \rfloor}]
& \le & \bbP[ H | A_{\lfloor d/2 \rfloor}] \\
 & \le  & \Big(1- \frac{C_9}{n^{({1/2}+
\eps)(d-2)}}\Big)^{r_{\lfloor d/2 \rfloor} n^{\lfloor d/2 \rfloor}}.
\end{eqnarray*}
 Now, choosing
$ \eps < \frac{1}{2(d-2)} $, and using the fact that
\begin{eqnarray*}
\lefteqn{ \bbP[ T(0,x_0) > n + (I+d^2)n^{1+2 \eps}] }\\
&\le &
\bbP[ T(0,x_0) > n + (I+d^2)n^{1+2 \eps} | A_{\lfloor d/2 \rfloor} ] +
\bbP[A_{\lfloor d/2 \rfloor}^c]
\end{eqnarray*}
together with~(\ref{condback}),
 we are finished. \qed

\medskip
\noindent
{\it Remark.} The sub exponential estimate for the tail of
the distribution of $ T(0,x) $ also holds for $ d = 1$.
The proof is similar to what is done for $ d = 2 $ and
therefore is omitted.

\section{Asymptotic shape}
\label{AS}

In the previous section we proved that for all  $ x \in \bbZ^d $,
the sequence $ (T(nx,(n+1)x), n \ge 0) $ satisfies the
hypotheses of Liggett's subadditive ergodic theorem. Therefore,
defining $ T(x) := T(0,x) $ for all $ x \in \bbZ^d $, it holds
that there exists $\mu(x) \geq 0$ such that
\begin{equation}
\label{limtint} \frac{T(nx)}{n} \to \mu(x) \qquad
    \mbox{a.s., $n\to\infty$.}
\end{equation}
   From the fact $ T(nx) \ge n\|x\|_1 $
it follows that $ \mu(x) \ge \|x\|_1 $ for
all $ x \in \bbZ^d$.

\begin{lem}
\label{largedeviation2}
For all $ x \in \bbZ^d $ there are constants $ 0 < \delta_0 < 1 $,
$ C > 0 $ and $ 0 < \gamma < 1 $ such that
\[
\bbP\Big[T(x) \geq \frac{\|x\|_1}{\delta_0}\Big] \le C \exp\{-
\|x\|_1^{\gamma}\} .
\]
\end{lem}

\noindent {\it Proof.} Let $ n := \|x\|_1 $ and $ 0 = x_ 0 , x_1 ,
x_2 , \ldots, x_n = x $ be a path connecting the origin to site $ x $
such that for all $ i $, $ \|x_i - x_{i-1}\|_1 = 1 $. Let $ Y_i :=
T(x_{i-1},x_i)$. Due to the subadditivity, it is enough to proof that
\begin{equation}
\label{sumyi}
 \bbP\Big[ \sum_{i=1}^n Y_i \geq \frac{\|x\|_1}{\delta_0}  \Big]  \le
C\exp\{-\|x\|_1^{\gamma}\}.
\end{equation}
Let
\[
 B := \Big\{ Y_i < \frac{\sqrt n}{2}, i=1, \ldots, n \Big\} .
\]
Clearly, by Theorem~\ref{expdecay} we have
\begin{equation}
\label{gshtrih}
 \bbP[ B ] \ge 1 -
C_1 n \exp\{-n^{\gamma'}\}
\end{equation}
For some $\gamma'>0$.
For $i=1,\ldots, \lceil {\sqrt n} \rceil$ let
\[
 \sigma_i := \sum_{j=0}^{M_i} Y_{i+j{\lceil {\sqrt n} \rceil}},
\]
where
\[
 M_i := \max \{ j \in \bbN :  i + j \lceil
{\sqrt n} \rceil \le n \} .
\]
 Observe that, if the event $ B $ happens, then each
$ \sigma_i $ is as a sum of independent
identically distributed
random variables, since in this situation the variables
$ \{ Y_{i+j{\lceil {\sqrt n} \rceil}}: j=1, \ldots, M_i \} $
depend on disjoint sets of random walks.

We point out that we cannot guarantee the existence of the moment
generating function of $ Y_i$. All we have is a
sub exponential estimate as in~(\ref{subexponential})
(see Theorem~\ref{expdecay}).
Lemma~\ref{Nagaev} takes care of the situation and allows us to
obtain (for $ \delta_0 = {1 / a} $)
\begin{eqnarray}
\bbP\Big[ \sum_{i=1}^n Y_i > \frac{n}{\delta_0} \Big| B \Big]
&\le& \bbP\Big[\Big\{\sigma_1 \le \frac{M_1}{\delta_0}, \ldots,
\sigma_{\lceil {\sqrt n} \rceil} \le
\frac{M_{\lceil {\sqrt n} \rceil}}{\delta_0}\Big\}^c \Big| B \Big]
\nonumber \\
&\le&  \sum_{i=1}^{\lceil {\sqrt n} \rceil} \bbP\Big[\sigma_i
  \geq \frac{M_i}{\delta_0} \Big| B \Big]
\nonumber \\
&\le& C_2 \sum_{i=1}^{\lceil {\sqrt n} \rceil}
\exp\{-M_i^{\beta}\}. \label{uravn4.3}
\end{eqnarray}
Note that if $ n $ is large then for all
$ i \le \lceil {\sqrt n} \rceil $
it holds that $ M_i = {\cal O}({\sqrt n}) $.
The result follows from~(\ref{gshtrih}) and~(\ref{uravn4.3}). \qed

\medskip

Let us extend the definition of $ T(x,y) $ to the
whole $ \bbR^d \times \bbR^d $ by defining
\[
 {\bar T}(x,y) = \min\{n:y \in {\bar \xi}^{x_0}_n\},
\]
where $x_0\in \bbZ^d$ is such that $x\in (-1/2,1/2]^d + x_0$.
Note that the subadditive property holds for $ {\bar T}(x,y)$ as
well. The next goal is to show that $ \mu $ can be extended to
$ \bbR^d $ in such a way that~(\ref{limtint}) holds for all $ x \in
\bbR^d $. As we did before, let us consider $ {\bar T}(x) = {\bar
T}(0,x)$.

\begin{lem}
\label{l4.2}
For all $ x \in \bbQ^d $
\[
 \frac{{\bar T}(nx)}{n} \to \frac{\mu(mx)}{m} ,
\]
where $ m $ is the smallest positive integer such that
$ mx \in \bbZ^d $.
\end{lem}

\noindent {\it Proof.} Let $ n = km + r $, where $ k, r \in \bbN $
and $ 0 \le r < m $. Since $ {\bar T}(nx) \le T(kmx) + {\bar T}(rx)$,
it is true that
\begin{equation}
\label{ura*}
 \limsup_{n \to \infty} \frac{{\bar T}(nx)}{n} \le \frac{\mu(mx)}{m}.
\end{equation}
Analogously, writing $n=(k+1)m-l$ one gets
$ T((k+1)mx) - {\bar T}(lx) \le {\bar T}(nx)$,
which implies that
\begin{equation}
\label{ura**}
 \liminf_{n \to \infty} \frac{{\bar T}(nx)}{n} \ge \frac{\mu(mx)}{m} .
\end{equation}
Combining~(\ref{ura*}) and~(\ref{ura**}), we finish
the proof of Lemma~\ref{l4.2}.
\qed

\medskip

By standard methods one can prove that $ \mu $ is uniformly continuous
in $ \bbQ^d $, and therefore can be continuously extended to $ \bbR^d $
in such a way that
\begin{equation}
\label{addit}
 \lim_{n \to \infty} \frac{{\bar T}(nx)}{n} = \mu(x).
\end{equation}
So, it follows that $\mu$ is a norm in $ \bbR^d $.

For $y\in\bbR^d$ and $a>0$ denote $D(y,a) = \{x\in\bbR^d :
   \|x-y\|_1 \leq a\}$.
\begin{lem}
\label{smallball}
There exist $ 0 < \delta < 1 $, $C>0$, $\gamma_0>0$ such that
\[
 \bbP[ D(0, {n \delta})
\subset {\bar \xi}_n ] \geq 1-C\exp\{-n^{\gamma_0}\}
\]
for all~$n$ large enough.
\end{lem}

\noindent {\it Proof.} By Lemma~\ref{largedeviation2},
if~$x$ is such that $\|x\|_1=n$ then we have
\[
  \bbP\Big[T(x) \geq \frac{n}{\delta_0}\Big]
 \leq C \exp\{-n^{\gamma}\}.
\]
Now, let $y\in \bbZ^d$ be such that $\|y\|_1<n$.
 Then there exists~$x\in \bbZ^d$
such that $\|x\|_1=n$ and $\|x-y\|_1=n$. As $T(y)\leq T(x)+T(x,y)$
and by Lemma~\ref{largedeviation2} we have
\begin{eqnarray*}
\bbP\Big[T(y) \geq \frac{2n}{\delta_0}\Big] &\leq &
 \bbP\Big[T(x)+T(x,y) \geq \frac{2n}{\delta_0}\Big]\\
&\leq & \bbP\Big[T(x) \geq \frac{n}{\delta_0}\Big] +
\bbP\Big[T(x-y) \geq \frac{n}{\delta_0}\Big]\\
&\leq & 2C\exp\{-n^{\gamma}\}.
\end{eqnarray*}
So,
\[
\bbP\Big[\mbox{there exists }y\in D(0,n)\mbox{ such that }
 T(y) \geq \frac{2n}{\delta_0} \Big] \leq C_1 n^d \exp\{-n^{\gamma}\},
\]
which finishes the proof of Lemma~\ref{smallball}.
\qed

\medskip

Now we are able to finish the proof of
     the shape theorem for the frog model.

\medskip

\noindent {\it Proof of Theorem~\ref{maintheo}}. \
Let $ { \mathsf A } := \{x \in \bbR^d : \mu(x) \le 1 \} $. We first
prove that
\[
 n(1-\eps){\mathsf A} \subset {\bar \xi}_n \ \hbox{ for all } n \hbox{ large enough,
 almost surely.}
\]

Since ${\mathsf A}$ is compact, there exist
$ {\mathsf F}:= \{x_1, \ldots, x_k\} \in {\mathsf A} $ such that
$ \mu(x_i) < 1 $ for $i=1,\ldots,k$, and
(with~$\delta_0$ from Lemma~\ref{largedeviation2})
\[
 {\mathsf A} \subset \bigcup_{i=1}^{k}D (x_i, \eps\delta_0).
\]
Note that~(\ref{addit}) implies that
$ n{\mathsf F} \subset
{\bar \xi}_{n}$ for all $n$ large
enough almost surely.

By Lemma~\ref{smallball}, we have almost surely that there exists
$ n_o $ such that for all $ n \ge n_0 $
\[
n D( (1-\eps) x_i, \eps\delta_0) \subset
{\bar \xi}^{ n (1-\eps) x_i}_{n \eps},
\hbox{ for all } i=1,2, \dots k
\]
 and this part of the proof is done.

Now we prove that
\[{\bar \xi}_n
\subset n(1+ \eps){\mathsf A} \
\hbox{ for all $n$ large enough almost surely}.
\]

First choose $ {\mathsf G} := \{y_1, \ldots, y_k\} \subset
2{\mathsf A}\setminus {\mathsf A} $
such that
\[
  2{\mathsf A}\setminus {\mathsf A}
\subset \bigcup_{i=1}^k D(y_i, \eps(1+\eps)^{-1}\delta_0 ) .
\]
Notice that $ \mu(y_i) > 1 $ for $i=1,\ldots,k$.
Analogously, (\ref{addit}) implies that

\begin{equation}
\label{um}
  n {\mathsf G} \cap {\bar \xi}_{n} =\emptyset \
\hbox{ for all $n$ large enough almost surely}.
\end{equation}

Suppose that, with positive probability,

\begin{equation}
\label{dois}
{\bar \xi}_{n} \not\subset
n(1+\eps){\mathsf A} \ \hbox{ for infinitely many } n \in \bbN
\end{equation}

\noindent
Fixed a realization of the process such that~(\ref{um})~and~(\ref{dois})
happens, choose $n_0$ so large that for $n>n_0,$ $ n {\mathsf G}
\cap {\bar \xi}_{n} =\emptyset$ and choose $n>n_0$ such that ${\bar \xi}_{n}
\not\subset n(1+\eps){\mathsf A}.$

Since  ${\bar \xi}_{n}$ is connected and~(\ref{dois}) holds, there is a site
$x \in {\bar \xi}_n \cap(1+\eps)n(2{\mathsf A}\setminus{\mathsf A}).$
By Lemma~\ref{smallball},
we can suppose that in the realization we are considering, $n_0$ is so large
that for $n>n_0,$  ${\mathsf D}(x,n\eps\delta_0) \subset {\bar \xi}_{n(1+\eps)}.$

Notice that, since  $(1+\eps)n(2{\mathsf A}\setminus {\mathsf A})
\subset \bigcup_{i=1}^k D((1+\eps)ny_i,n \eps\delta_0 ),$ we must have

$$n(1+\eps)y_k \in {\mathsf D}(x,n\eps\delta_0) \subset {\bar
\xi}_{n(1+\eps)}.$$

\noindent
This contradicts~(\ref{um}), and, therefore, concludes the proof of the theorem. \qed

\medskip

\noindent {\it Proof of Theorem~\ref{flat_edge}.} To prove the
theorem, it is enough to prove the following fact: for
fixed~$i,j$, there exists~$\beta$ such that
\begin{equation}
\label{betaij}
 \Theta_{ij}^\beta \subset {\sf A}_m.
\end{equation}
Indeed, in this case~(\ref{betaij}) holds for all~$i,j$ with the
same~$\beta$ by symmetry, hence $\Theta^\beta \subset {\sf A}_m$
by virtue of convexity of~${\sf A}_m$.

Now, the proof of~(\ref{betaij}) is just a straightforward
adaptation of the proof of ``flat edge'' result of~\cite{DL}. To
keep the paper self-contained, let us outline the ideas of the
proof. Suppose, without loss of generality, that $i=1, j=2$
and~$m$ is even. We are going to prove that the frog model
observed only on $\Lambda_{12}\cap \bbZ^d_+$ dominates the
oriented percolation process in~$\bbZ^2_+$ with parameter
$\theta=1-(1-(2d)^{-1})^{m/2}$. To show this, first suppose that
initially for any~$x$ all the particles in~$x$ are labeled
``$x{\to}$'' or ``$x{\uparrow}$'' in such a way that~$x$ contains
exactly~$m/2$ particles of each label. Define $e_1,e_2$ to be the
first two coordinate vectors. The oriented percolation is then
defined in the following way: For $x\in\Lambda_{12}\cap \bbZ^d_+$
\begin{itemize}
\item the bond from~$x$ to~$x+e_1$ is open if for the frog model
at the moment next to that of activation of the site~$x$ at least
one particle labeled ``$x{\to}$'' goes to~$x+e_1$;
\item the bond from~$x$ to~$x+e_2$ is open if at that moment at least one
particle labeled ``$x{\uparrow}$'' goes to~$x+e_2$.
\end{itemize}
Clearly, the two above events are independent, and their
probabilities are exactly~$\theta$. So the frog model indeed
dominates the oriented percolation in the following sense: if a
site $x=(x^{(1)},x^{(2)})$ (for the sake of brevity forget the
zero coordinates from~$3$ to~$d$) belongs to cluster of~$0$ in the
oriented percolation, then in the frog model the corresponding
site is awakened {\it exactly\/} at time $x^{(1)}+x^{(2)}$. Now it
rests only to choose~$m$ as large as necessary to make the
oriented percolation supercritical ($\theta\to 1$ as $m\to
\infty$) and use the result that (conditioned on the event that
the cluster of~$0$ is infinite) the intersection of the cluster
of~$0$ with the line $\{(x^{(1)},x^{(2)}):x^{(1)} + x^{(2)}=n\}$
grows linearly in~$n$ (cf.~\cite{DL}). \qed

\medskip
\noindent {\it Proof of Theorem~\ref{full_diamond}.} Denote
\[
{\cal D}_n = \{x \in \bbZ^d : \|x\|_1 \leq n\}.
\]
Choose $\theta<1$ such that $\delta<\theta d$. Start the process
and wait until the moment~$n^\theta$. As~$\eta\geq 1$ a.s., by
Lemma~\ref{smallball} there exists~$C_1$
such that with probability at least
$1-\exp\{n^{\gamma_0}\}$ all the frogs in the ball of
radius~$C_1n^\theta$ centered in~$0$ will be awake. Let
$x_1,\ldots,x_N$ be the sites belonging to that ball enumerated in
some order. Let~$\zeta_i$ be the indicator of the following event:
\[
 \{\mbox{in the initial configuration $x_i$ contains less than
$(2d)^n$ particles}\}.
\]
Clearly, the inequality~(\ref{fd}) implies that $\bbP[\zeta_i=1]
\leq 1-n^{-\delta}$. As~$N$ is of order~$n^{\theta d}$, one can
apply Lemma~\ref{Shiryaev} with $p=1-n^{-\delta}$ and
$a=1-n^{-\delta}/2$ to get that with probability at least
$1-\exp\{-C_2n^{\theta d-\delta}\}$ at time~$n^\theta$ one will
have at least $C_3 n^{\theta d-\delta}$ activated sites~$x_i$
with~$\zeta_i=0$. This in turn means that
 at time~$n^\theta$ there are at least $C_3 n^{\theta d-\delta} (2d)^n$
active frogs in~${\cal D}_{n^\theta}$. Note the following simple fact:
If~$x$ contains at least $(2d)^n$ active particles and
$\|x-y\|_1\leq n$,
then until time~$n$ with probability bounded away from~$0$ at
least one of those particles will hit~$y$. Using this fact, we get
that with overwhelming probability all the frogs in the diamond
${\cal D}_{n-n^\theta}$ will be awake at time $n^\theta+n$, which
completes the proof of Theorem~\ref{full_diamond}. \qed

\section{Remarks about continuous time}
A continuous-time version of the frog model can also
be considered. The difference from the discrete-time
frog model is of course that here the particles,
after being activated, perform a continuous-time
SRWs with jump rate~$1$. In the continuous-time context
Theorem~\ref{maintheo} also holds and its proof
can be obtained just by following the steps of our
proof for the discrete case. The only difficulty
that arises is that for continuous time, it is not so
evident that $\mu(x)$ (defined
by~(\ref{limtint})) is strictly
positive for $x\neq 0$, i.e.\ we must rule out the
possibility that the continuous-time frog model
grows faster than linearly. To overcome that difficulty,
note the following fact (compare with Lemma~9 on
page~16 of Chapter~1 of~\cite{D}): there exist a
positive number~$C$ such that, being $\|x\|_1\geq Cn$,
$\bbP[T(0,x)<n]$ is exponentially small in~$n$.
This fact by its turn easily follows from a domination of the frog model by
branching random walk.

\section*{Acknowledgements}
The authors have benefited on useful comments
and suggestion
from K.~Ravishankar and L.~Fontes
and wish to thank them.
E.~Kira should be thanked for a careful
reading of the first draft of this paper.
F.P.M. is also indebted to R.~Durrett
and M.~Bramson for
discussions about a continuous-time version of the model.

\end{document}